\documentclass[autocontact,english]{gaceta}

\usepackage[T1]{fontenc}
\usepackage{lmodern}

\usepackage[utf8]{inputenc}
\usepackage[spanish]{babel}

\usepackage{amsmath, amsthm, amssymb}
\usepackage{color}

\journame{La Gaceta de la RSME}
\yearofpublication{2020}
\volume{23}
\issuenumber{2}
\modifiedPaper
\modifiedPaperInfo{Extended English version of the paper /
Versi\'on extendida en ingl\'es del art\'{\i}culo}
\originalPaperPages{243--261}
\belongstopart{Art\'{\i}culos}

\title{Remembering Louis Nirenberg and his mathematics}
\author{\large Juan Luis Vázquez, \\[4pt] \normalsize Real Academia de Ciencias, Spain} 
\shorttitle{Louis Nirenberg}

\contact{Juan Luis Vázquez, Departamento de Matemáticas,
Universidad Autónoma de Madrid \\ and Real Academia Española de Ciencias}
{juanluis.vazquez@uam.es}{http://verso.mat.uam.es/~juanluis.vazquez/}

\newcommand{\R}{\mathbb{R}}
\begin{document}

\selectlanguage{english}
\renewcommand{\porname}{by}
\renewcommand{\emailname}{Email}
\renewcommand{\webname}{Web page}
\maketitle

\begin{abstract}
The article is dedicated to recalling the life and mathematics of Louis Nirenberg, a distinguished Canadian mathematician who recently died in New York, where he lived. An emblematic figure of analysis and partial differential equations in the last century, he was awarded the Abel Prize in 2015. From his watchtower at the Courant Institute in New York, he was for many years a global teacher and master. He was a good friend of Spain.
\end{abstract}

\begin{center}
\includegraphics[width=0.8\textwidth]{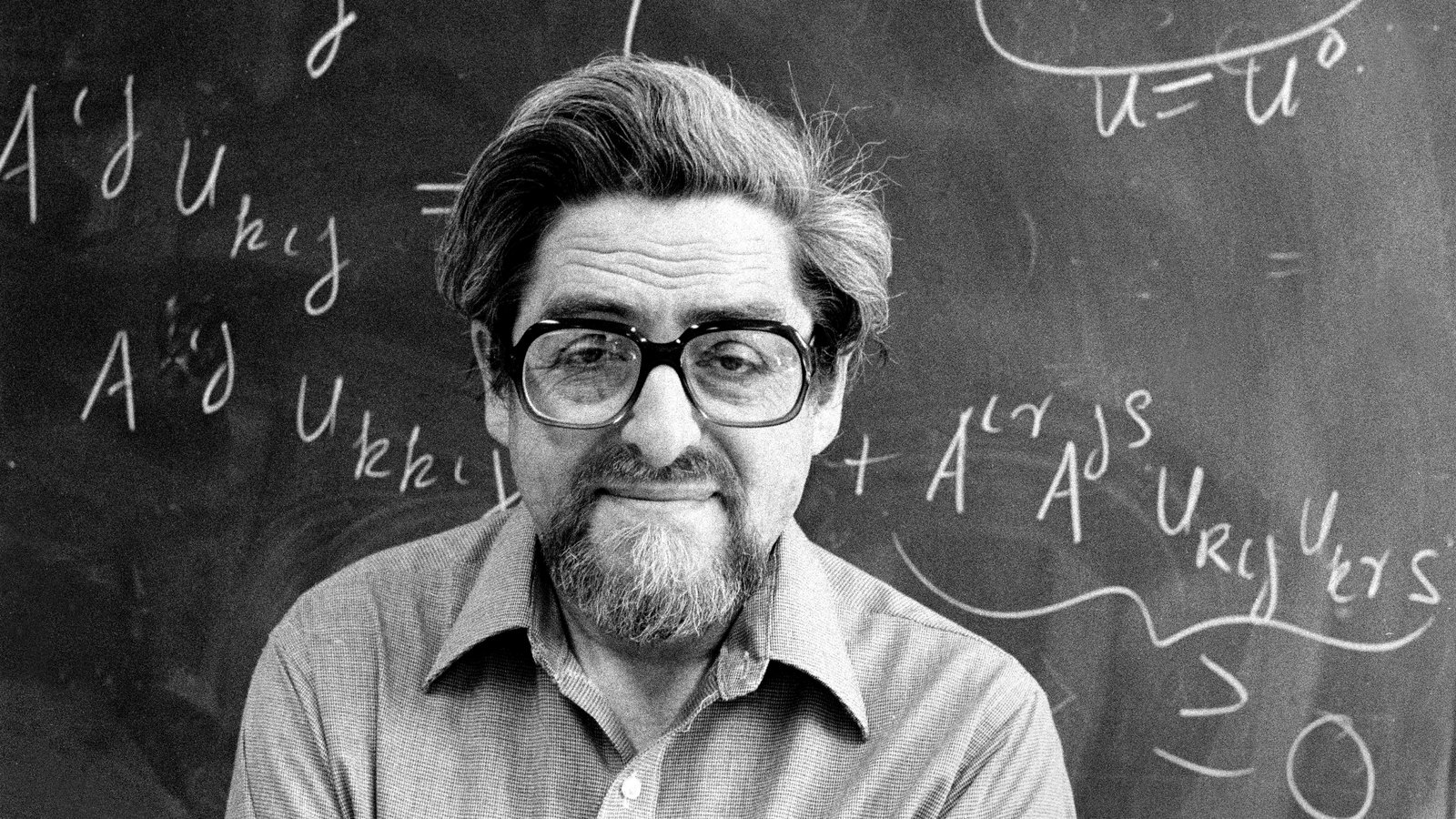}
\end{center}
\medskip
\hfill\begin{minipage}{0.75\textwidth}
\textit{One of the wonders of mathematics is you go somewhere in the world and you
meet other mathematicians, and it is like one big family. This large family is a
wonderful joy.}\footnotemark
\end{minipage}\qquad\quad\par
\footnotetext{From an interview with Louis Nirenberg appeared in  \textit{Notices of the AMS},
2002, \cite{AJack}}

\section{Introduction}

This article is dedicated to remembering the life and work of the prestigious Canadian mathematician Louis Nirenberg, born in Hamilton, Ontario, in 1925, who died in New York on January 26, 2020, at the age of 94. Professor for much of his life at the mythical Courant Institute of New York University, he was considered one of the best mathematical analysts of the 20th century, a specialist in the analysis of partial differential equations (PDEs for short).

When the news of his death was received, it was a very sad moment for many mathematicians, but it was also the opportunity of reviewing an exemplary life and underlining some of its landmarks. His work unites diverse fields between what is considered Pure Mathematics and Applied Mathematics, and in particular he was a cult figure in the discipline of Partial Differential Equations, a key theory and tool in the mathematical formulation of many processes in science, in engineering, and in other branches of mathematics. His work is a prodigy of sharpness and logical perfection, and at the same time its applications span  today multiple scientific areas.

In recognition of his work, in 2015  he received the Abel Prize along with the another great mathematician, John Nash. The Abel Prize is one of the greatest awards in Mathematics, comparable to the Nobel prizes in other sciences. At that time, the Courant Institute, where he was for so many decades a renowned professor, published an article called {\sl Beautiful Minds}\footnote{\url{https://www.nyu.edu/about/news-publications/news/2015/march/beautiful-minds-courantsnirenberg-princetons-john-nash-win-abel-prize-in-mathematics-.html}}
 which is quite enjoyable reading.

He was a distinguished member of the AMS (American Mathematical Society). Throughout his life, he received many other honors and awards, such as the AMS Bôcher Memorial Prize (1959), the Jeffery-Williams Prize (1987), the Steele Prize for Lifetime Achievement (1994 and 2014), the National Medal of Science (1995), the inaugural Crafoord Prize from the Royal Swedish Academy (1982), and the first Chern medal at the 2010 International Congress of Mathematicians, awarded by the International Mathematical Union and the Chern Foundation. He was a plenary speaker at the International Congress of Mathematicians held in Stockholm in August 1962; the title of the conference was \textit{``Some Aspects of Linear and Nonlinear Partial Differential Equations''\rm}. In 1969 he was elected Member of the U.S. National Academy of Sciences.

It was not honors that concerned him most, but rather his profession and the mathematical community that surrounded him. In his long career at the Courant Institute he discovered many mathematical talents and collaborated in numerous relevant works with distinguished colleagues. A wise man in science and life, he was one of the most influential and beloved mathematicians of the last century,  and the current century too. His teaching extended first to the international centers that he loved to visit, and then to the entire world. Indeed,  we live at this height of time  in a world-wide scientific society whose close connection brings so many benefits to the pursuit of knowledge. Many of his articles are among the most cited in the world.\footnote {Topic 35, PDEs, from the mathematical database Mathscinet, includes 3 articles by L. Nirenberg among the 10 most cited ever.}

\begin{figure}
\centering
\includegraphics[width=0.85\textwidth]{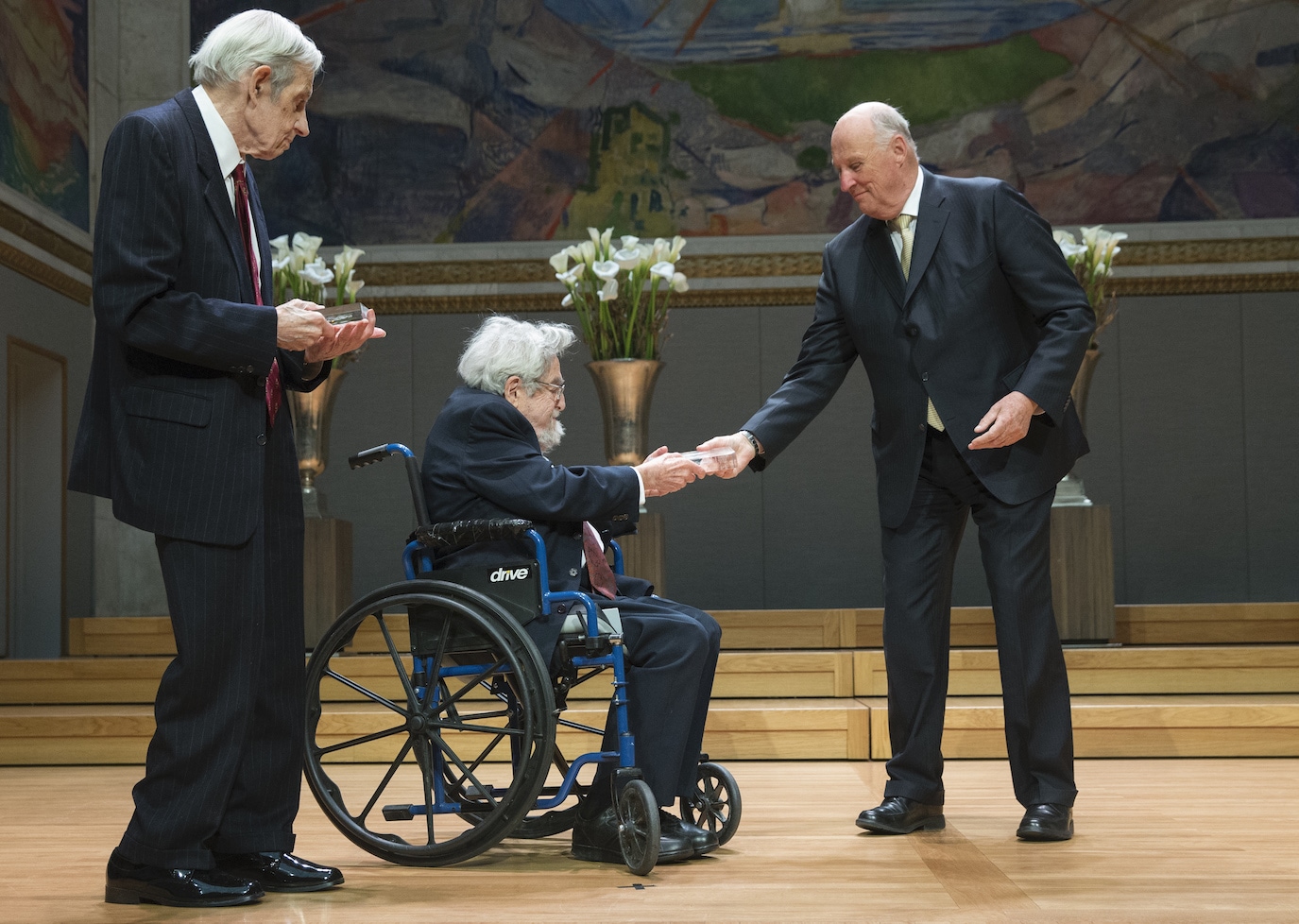}
\pseudocaption{Louis Nirenberg receiving the Abel Prize from King Harald V of Norway in the presence of John Nash (photo: Berit Roald/NTB scanpix).}
\end{figure}

\section{Starting}

In order to start the tour of his mathematics, nothing better than to quote a few paragraphs from the mention of the Abel Prize Committee in 2015:\footnote{See \url{https://www.abelprize.no/nyheter/vis.html?tid=63589}.}

\medskip

\noindent{\bf Mathematical giants: \sl Nash and Nirenberg are two mathematical giants of the twentieth century. They are being recognized for their contributions to the field of partial differential equations (PDEs), which are equations involving rates of change that originally arose to describe physical phenomena but, as they showed, are also helpful in analyzing abstract geometrical objects.\rm

The Abel committee writes: {\sl ``Their breakthroughs have developed into versatile and robust techniques that have become essential tools for the study of nonlinear partial differential equations. Their impact can be felt in all branches of the theory.''}

About Louis they say:   {\sl ``Nirenberg has had one of the longest and most f\^eted careers in mathematics, having produced important results right up until his 70s. Unlike Nash, who wrote papers alone, Nirenberg preferred to work in collaboration with others, with more than 90 per cent of his papers written jointly. Many results in the world of elliptic PDEs are named after him and his collaborators, such as the Gagliardo–Nirenberg inequalities, the John–Nirenberg inequality and the Kohn–Nirenberg theory of pseudo-differential operators.}

They conclude: {\sl ``Far from being confined to the solutions of the problems for which they were devised, the results proven by Nash and Nirenberg have become very useful tools and have found tremendous applications in further contexts''.}

To be precise, Nirenberg made fundamental contributions to  both linear and nonlinear partial differential equations, functional analysis, and their application in geometry and complex analysis. Among the most famous contributions we will discuss are the Gagliardo-Nirenberg interpolation inequality, which is important in solving the elliptic partial differential equations that arise within many areas of mathematics; the formalization of the BMO spaces of bounded mean oscillation, and others that we will be seeing.

A work of utmost relevance was the work with Luis Caffarelli and Robert Kohn aimed at solving the big open problem of existence and smoothness of the solutions of the Navier-Stokes system of fluid mechanics. This work was described by the AMS in 2002 as ``one of the best ever done.'' The problem is on the Millennium Problems List (the list compiled by the Clay Foundation), and is one of the most appealing open problems of mathematical physics, raised nearly two centuries ago. Fermat's Last Theorem and the Poincaré Conjecture have been defeated at the turn of the century, but the Navier-Stokes enigma (and in some sense its companion about the Euler's system) keep defying us. We will deal with the issue in detail in Section \ref{sec.ns}.


\subsection{The beginnings. From Canada to New York}

Louis Nirenberg grew up in Montr\'eal, where his father was a Hebrew teacher. After graduating\footnote{With a degree in mathematics and physics, also in mathematics being bilingual counts.} in 1945 at McGill University, Montréal, Louis found a summer job at the National Research Council of Canada, where he met the physicist Ernest Courant, the son of Richard Courant, a famous professor at New York University. Ernest mentioned to Nirenberg that he was going to New York to see his father and Louis begged him for advice on a good place to apply for a master study in physics. He returned with Richard Courant's invitation for Louis to go to New York University (NYU) for a master's degree in mathematics, after which he would be prepared for a physics program.

But once Louis began studying Mathematics at NYU, he never changed. He defended his doctoral thesis under James Stoker in 1949, solving a problem in differential geometry. The dice were cast. We  reach a crucial moment in Louis's life. Breaking with the golden rule\footnote{which is an essential part of the American professional practice.} according to which ``a recent doctor should move to a different environment'', Richard Courant kept his best students around him, including Louis Nirenberg, and he thus created the NYU Mathematical Institute, the famous Courant Institute, which has become a world benchmark for high mathematics, comparable only to the Princeton Institute for Advanced Study on the East Coast of the USA. Louis was first a postdoc and then a permanent member of the faculty. There he thrived and spent his life.

\subsection{Equations and Geometry}

The problem Stoker gave to Louis   for his thesis, entitled {\sl ``The Determination of a Closed Convex Surface Having Given Line Elements''}, is called ``the embedding problem'' or ``Weyl Problem''. It can be  stated as follows: Given a 2-dimensional sphere with a Riemannian metric such that the Gaussian curvature is positive everywhere, the question is whether a surface can be constructed in three-dimensional space so that the Riemannian distance function coincides with the distance inherited from the usual Euclidean distance in the three dimensional space (in other words, whether there is an isometric
embedding as a convex surface in $ \R^3 $). The great German mathematician Hermann Weyl had taken a significant first step to solve the problem in 1916, and Nirenberg, as a student, completed Weyl's construction. The work to do was to solve a system of nonlinear partial differential equations of the so-called ``elliptical type''. It is the kind of equation and application that Louis Nirenberg has been working on ever since. Progress has been slow but continued over time and is impressive at this moment.\footnote{Isometrically embedding low dimensional manifolds into higher dimensional Euclidean spaces is the contents of a famous paper by J. Nash in 1956.}

\section{The power and beauty of inequalities}

Let us focus our attention on one of the most relevant topics in Louis Nirenberg's broad legacy, at the same time closest to our mathematical interests. (Almost) every career in PDEs begins with the study of linear elliptic equations. These form nowadays a well-established theory  which combines Functional Analysis, Calculus of Variations, and explicit representations to produce solutions in suitable functional spaces. For the classical equilibrium equations in the mechanics of continuous media, known as Laplace's and Poisson's equations, in symbols $ \Delta u = f $, there is a classic ``maximum principle'' that provides the necessary estimates that guarantee the existence and uniqueness of solutions. When combined with skillful tricks of the trade, it makes possible to obtain finer estimates, such as regularity and other properties. Let us mention the estimates known under the names Harnack and Schauder,   cf. \cite{Evansbk, GT77}. In this regard, Nirenberg is quoted as saying, either jokingly or seriously,

\smallskip

\centerline{\sl ``I made a living off the maximum principle''.\footnote{Curiously, it applies to me too. My most read article deals with the ``Strong Maximum Principle'', \cite{VSMP}.}}

\smallskip


Many of the interesting problems that are proposed in Physics and other sciences and involve PDEs  are {\bf nonlinear,} for example the fluid equations or the curvature problems in geometry. These nonlinear problems can seldom be solved by explicit formulas. Because of that difficulty, the mathematical study of these problems has attracted increasing attention from the best mathematical minds of the past century, with remarkable success stories. The usual approach goes as follows: the solution has to be obtained by some kind of approximation, and an essential technical point is usually to show that the proposed approximation procedure (or procedures) converge to a solution\footnote{Taken in some sense acceptable to physics, for example, the solution in the weak sense or the solution in the distributional sense.}. A complicated topology and functional analysis machinery has been developed over time and is available to test such convergence, provided certain estimates are fulfilled; thir role is to allow for the approximation to be controlled. See in this sense the book that many of us have studied as young people \cite{BrzAF}.

Much of the work of an ``EPD Analyst''\footnote{\sl Analysis of PDEs \rm is an area of Mathematics in the US that perfectly describes our specialty  which is neither pure nor applied, and does not need to declare itself in either direction. Such a denomination is not much used in Spain and other countries;  that is, in my opinion, the source of some persistent malfunctions.}
consists in finding estimates that control the passage to the limit that has to be applied, or to find a convenient fixed point theorem. A common saying in our trade goes as follows: {\sl Existence theorems come from a priori estimates and suitable functional analysis}.  Estimate, this is the key word in the world that Louis Nirenberg and his colleagues bequeathed us. ``Estimate'' means the same thing as ``inequality'', and here we refer of course to a functional or numerical inequality.

It may look surprising  to the reader, even weird, to find it so clearly stated: Inequalities, and not equalities (or identities), are the technical core of such a central theory of mathematics as PDEs. However,  this is precisely the mathematical revolution that was in the making when Louis was young. Indeed, when he arrived at NYU, the most active and renowned researcher was probably Kurt Otto Friedrichs, who decisively influenced Nirenberg's future research career. Friedrichs loved inequalities, as Louis put it:
\begin{center} {\sl ``Friedrichs was a great lover of inequalities and that affected me \\ very much. The point of view was that the inequalities are \\ more interesting than the equalities'' \rm.}
\end {center}
Carrying forward on that idea, Nirenberg has been unanimously recognized as a ``world master of inequalities''. Here is another saying by Louis:
\begin{center}
\sl ``I love inequalities. So if somebody shows me a new inequality, I say: \\ "Oh, that's beautiful, let me think about it," and \\I may have some ideas connected to it''.
\end{center}

For many years, mathematicians from all over the world came to the Courant Institute to seek his advice on issues involving inequalities.

And there we are. We do not reject or despise the beauty of the exact solution if there is one, but functional inequalities are our firm support in an uncertain world that is yet to be discovered and described. The key technical point of modern PDE theory is to establish the most needed and appropriate estimates in the strongest possible way.

\subsection{Sobolev, Gagliardo and Nirenberg}

There are many types of estimates the researcher needs in the study of nonlinear PDEs, but some have turned out to be much more relevant than others. We will talk here about a type that has become particularly famous and useful. They are often collectively called ``Sobolev estimates'' in honor of the great Russian mathematician Sergei Lvovich Sobolev because of his seminal work \cite{Sob38}, 1938. Briefly stated, they  estimate the norms of functions belonging to the Lebesgue spaces $ L^p (\Omega)$, $ 1 \le p \le \infty $, in terms of their (weak) derivatives of various orders.  In 1959 Emilio Gagliardo \cite{Gagl59} and Louis Nirenberg \cite{Nir59} gave an independent and very simple proof of the following inequality:

\medskip

\noindent {\bf Theorem \rm  (Gagliardo-Nirenberg-Sobolev Inequality)}. {\sl Let $1\le p < n$. There exists a constant $C>0$ such that the following inequality
$$
\|u\|_{L^{p*}(\R^n)}\le C\|Du\|_{L^p(\R^n)},\quad p^*:=np/(n-p),
$$
holds true for all functions  $u\in C^1_c(\R^n)$. The constant $C$ depends only on $p$ and $n$. The exponent $p^*$ is called the Sobolev conjugate of $ p$. $Du$ denotes the gradient vector.}

\medskip

Gagliardo and Nirenberg included as their starting point the important case of exponent $ p = 1$, left out by Sobolev. The inequality implies the continuous inclusion of the Banach space called $W^{1,p}(\R^n)$ into $L^{p^*}(\R^n)$ (immersion theorem). Versions for functions defined  in bounded open sets $\R^n$ followed naturally. This inequality soon attracted multiple applications and a wide array of variants and improvements. Very interesting versions deal with functions  defined on Riemannian manifolds. We comment below  four  additional aspects that we find appropriate for the curious reader.

\begin{itemize}
\item[(i)] Thierry Aubin \cite{Aub76} and Giorgio Talenti \cite{Tal76} obtained in 1976 the \sl best constant \rm in this inequality, finding the functions that exhibit the \sl worst behaviour\rm\footnote{This is an apparent grammatical contradiction that gives rise to beautiful functions.}. Indeed, when $ 1 <p <n $ the maximum quotient $\|u\|_{L^{p*}(\R^n)}/\|Du\|_{L^{p}(\R^n)} $ is optimally realized by the function
$$
U(x)= \left(a+ b\, |x|^{p/(p-1)}\right)^{-(n-p)/p}
$$
where $ a, b> 0 $ are arbitrary constants.\footnote{We ask the reader to consider the simple case $ a = b = 1 $, $ p = 2 $ in dimension $ n = 4$. The function looks a bit like Gaussian but it is not at all.} It is the famous \sl Talenti profile\rm. Note that $(n-p)/p=n/p^*$. It happens that $ U $ is a probability density (integrable) if $ (n-p) / (p-1)> n $, that is, if $ 1 <p <p_c = 2n / (n + 1) $. The $ U $ profile and its powers appear recurrently in PDEs. Thus, in nonlinear diffusion we find it as a power of the Barenblatt profile in fast diffusion, see Chapter 11 of \cite{VazSmooth}, and the curiously critical exponent $ p_c $ also appears, but with consequences that go in the converse direction.

\begin{figure}
\centering
\includegraphics[width=0.4\textwidth]{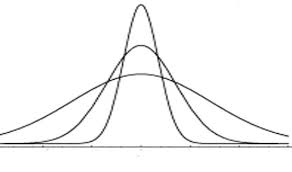}
\pseudocaption{The Talenti profile for different values of the parameters.}
\end{figure}

\item[(ii)] Gagliardo and Nirenberg's work extends to the famous {\sl Gagliardo-Nirenberg interpolation inequality}, a result in  Sobolev's theory of spaces that estimates a certain norm of a function in terms of a product of norms of functions and derivatives thereof. We enter here a realm of higher complexity.\footnote{We will avoid further details on these inequalities that can be found in the cited references.} See details in \cite{BM2018}.

\item[(iii)]  In 1984 Luis Caffarelli, Bob Kohn and Louis Nirenberg needed inequalities of the previous type in functional Lebesgue spaces but with the novelty of including so-called \sl weights\rm, and this motivated the article \cite{CKN84} on the famous {\sl CKN estimates originated} for spaces with power weights. This was the beginning of an extensive literature. A very striking effect arose in those studies: unlike GNS inequalities, there exists a phenomenon of symmetry breaking in the CKN inequalities, i.\,e., minimizers of such inequalities need not be symmetric functions, even when posed in the whole space or in balls. The exact range of parameters for the symmetry breaking was found by J. Dolbeault, M. J. Esteban and M. Loss in~\cite{DEL-2016}.

\item[(iv)]  In 2004 D. Cordero-Erausquin, B. Nazareth and C. Villani \cite{CENV04} used mass transport methods to obtain sharp versions of the Sobolev-Gagliardo-Nirenberg inequalities. Mass transport is one of the most powerful new instruments used in PDE research. This topic is related to the  isoperimetric inequalities of ancient fame,\footnote {See
\tt https://en.wikipedia.org/wiki/Isoperimetric-inequality.} that  now live moments of fruitful coincidence with  Sobolev theory. The survey \cite{XC17} talks about this relationship.
\end{itemize}

The world of estimates that we have outlined has came to be an enormous space presided over by quite distinguished  names, like  H. Poincaré, J. Nash, G. H. Hardy, C. Morrey, J. Moser, N. Trudinger and other remarkable figures.  Hardy-Littlewood-Pólya's book \cite{HardyLP} had a great influence on generations of analysts. A commendable book on the importance of inequalities in Physics is the second volume of Elliott Lieb's selected works, \cite{Lieb}.

As a representative  example chosen from among the numerous recent works, I would like to mention the article by M. del Pino and Jean Dolbeault \cite{DPDo02}. It  establishes a new optimal version of the Euclidean  Gagliardo-Nirenberg inequalities. This allows the authors to obtain the convergence rates to the equilibrium profiles of some nonlinear diffusion equations, such as those of the ``porous media'' type, one of the leitmotifs of my research. The authors completed the study and application with two new articles in 2003. New functional inequalities based on entropy, maximum principles, and symmetrization processes allowed a group of us to find convergence rates for very fast diffusion equations in \cite{BBDGV}, thus solving in 2009 a much studied open problem. It was almost 3 years of work by a team of 5 people. Plus the work of  previous authors.

Finally, there is  a great deal of activity in the world of Sobolev spaces of fractional order (also called Slobodeckii spaces),   and the associated fractional diffusions, cf. \cite{Hitch, CPFMV}. It is a topic in full swing, a part of my current mathematical efforts.\footnote{There is a wide representation of Spanish mathematicians active in these subjects with remarkable results that would be well worth a review.}

\subsection{New spaces. John-Nirenberg space}

Let us go back for a moment to the origins. The limiting case of the Gagliardo-Nirenberg-Sobolev inequality happens for $ p = n$. Thanks to new inequalities due to C. Morrey, we know that for $ p> n $ the resulting functions are H\"older continuous  functions, \cite{Evansbk}. But the $ p = n $ case was bizarre and it was left to Fritz John and Louis Nirenberg to solve the puzzle in 1961  by introducing the new BMO space of  functions of {\sl bounded mean oscillation}, see \cite{JN61}. Actually, BMO is not a function space but rather a space of function classes modulo constants. For this space there is the appropriate inequality.

\medskip

\noindent {\bf Theorem \rm (John-Nirenberg)}. {\sl If $u \in W^{1,n}(\R^n)$ then $u$  belongs to BMO and
 $$
 {\displaystyle \|u\|_{BMO}\leq C\|Du\|_{L^{n}(\mathbf {R} ^{n})},}
$$
for a constant $ C>0$ depending only on $ n.$}\footnote{The curious reader will wonder which function optimizes the constant. So?}

The BMO spaces are one a very popular new object in functional and harmonic analysis, they replace $ L^\infty $ when it turns out so. They were characterized by Charles Fefferman in \cite{Feff71}. The BMO spaces are slightly larger than $ L^\infty $. The possible inequality (and functional immersion) of John-Nirenberg type using $ L^\infty $ instead of BMO as image space may seem reasonable but it is false.\footnote{The reader is kindly asked to find an elementary counterexample.} We ought to be very careful then with the critical cases, that Louis treated with utmost attention. The John-Nirenberg spaces are used in analysis, in partial differential equations, in stochastic processes, and in multiple applications. The reader may use the references    \cite{Leoni} and \cite{Brasco2021} for some updates to recent work.


\begin{figure}
\centering
\includegraphics[width=0.85\textwidth]{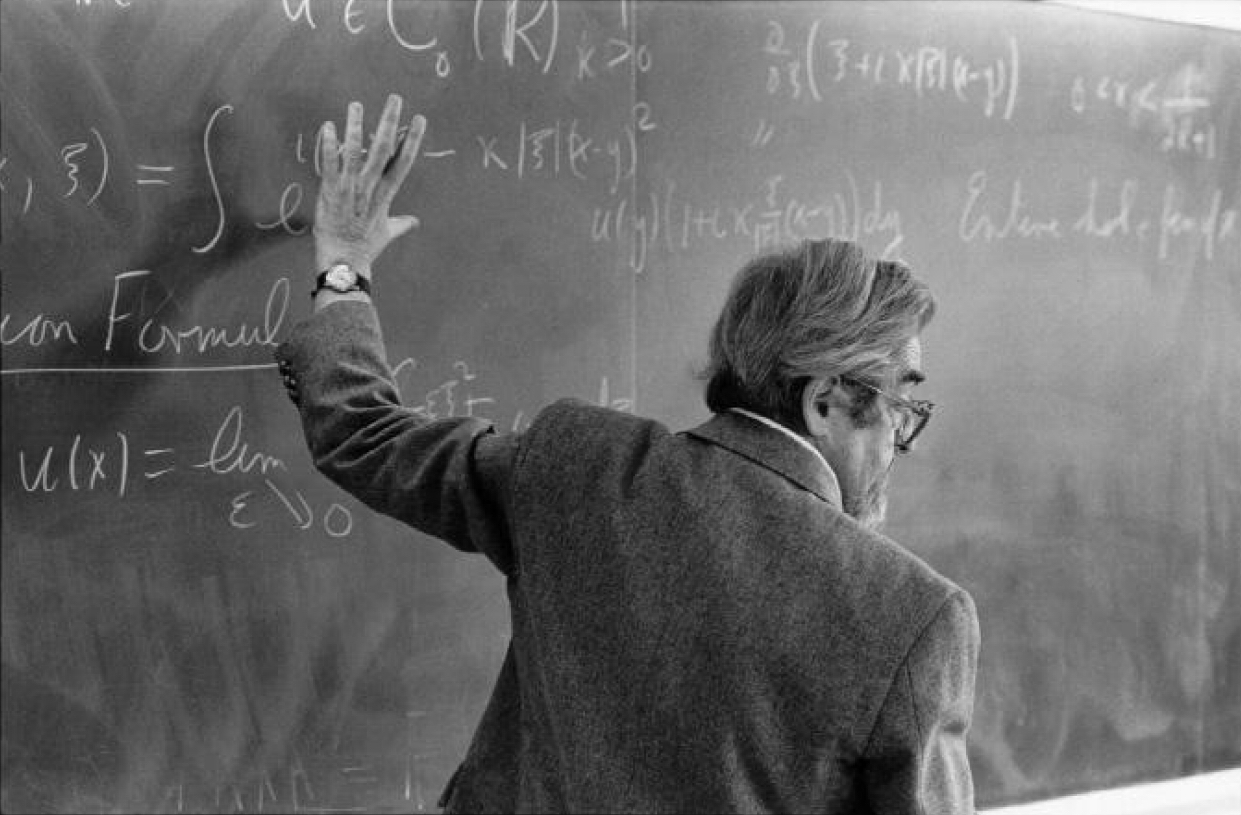}
\pseudocaption{Nirenberg on the blackboard (photo: Courant Institute, NYU).}
\end{figure}

\section{Navier-Stokes Equations}
\label{sec.ns}

The Navier-Stokes system of equations describes the dynamics of an incompressible viscous fluid. It was proposed in the 19th century to correct Euler's equations of ideal fluids, and adapt them to the more realistic viscous real world, \cite{Batch99}. The system reads
\begin{equation}
\label{e5.7}
\begin{gathered}
\dfrac{\partial \mathbf{u}}{\partial t} + (\mathbf{u}\cdot \nabla)\mathbf{u}
+ \dfrac1{\rho}\, \nabla p = \nu\,\Delta \mathbf{u} + \dfrac1{\rho}\,\mathbf{f}, \\
\nabla \cdot \mathbf{u}=0,
\end{gathered}
\end{equation}
where $\mathbf{u}$ is the velocity vector, $p$ is the pressure, both variable, while~$\rho$ (the density) and~$\nu$ (the viscosity) can be taken as positive constants. It has had a spectacular success in practical science and engineering, but its essential mathematical aspects (existence, uniqueness, and regularity) have offered a stubborn resistance in the physical case of three space dimensions (three or greater than three for the mathematician).

Fundamental works to cast the theory in a modern functional framework are due to Jean Leray \cite{L1, L2}, who already in 1934 speaks of weak derivatives in spaces of integrable functions.
Using the new methods of functional analysis, authors soon obtained estimates that proved to be good enough to establish the existence and uniqueness of Leray solutions in two space dimensions, $n = 2$. Furthermore, for regular initial data the solution is classical. But the advance stopped sharply in higher dimensions, $n \ge 3$.
We give the word to Charles Fefferman, of Princeton University, in his description of the open problem as the Clay Foundation Millennium Problem. It is about proving or refuting the following Conjecture:

\begin{quote}
\textit{(A) Existence and smoothness of Navier-Stokes solutions on $\R^3$. Take viscosity $\nu > 0$ and $n = 3$. Let $u_0 (x)$ be any smooth, divergence-free vector field satisfying the regularity and
decay conditions ({\rm specified}). Take external force $f(x,t)$ to be identically zero. Then there exist smooth functions $p(x,t)$, $u_i (x,t)$ on $\R^3 \times [0,\infty)$ that satisfy the Navier-Stokes system with initial conditions in the whole space.}\footnote{See full details of the presentation in \url{https://www.claymath.org/sites/default/files/navierstokes.pdf}.}
\end{quote}

The most significant advance in this field is in our opinion the article \cite{CKN82} in which L. Caffarelli, R. Kohn and L. Nirenberg attacked the problem of regularity and showed that if a solution with classical data develops singularities in a finite time, the set of such singularities must be in any case quite small in size. More specifically, ``the one-dimensional measure, in the Hausdorff sense, of the set of possible singularities (located in space-time) is zero''. This implies that if the singular set is not empty, it cannot contain any line or filament. In 1998 F. H. Lin \cite{Lin} gave an interesting new proof of this result.

We are talking about one of the milestones of the authors' career; it happened during the stay of a young Luis Caffarelli at the Courant Institute at Louis's invitation, and was published in 1982. The topic Fluids is completely different from the previous sections, but the functional estimates in Sobolev spaces play an essential role, along with the machinery of geometric measure theory.

The possible presence of these singularities was conjectured by Leray as a possible explanation for the phenomenon of \sl turbulence\rm. According to this hypothesis, even for regular data, solutions in three or more dimensions can develop singularities in finite time in the form of points where the so-called vorticity becomes infinite.

In the elapsed time, it has not been possible to prove or refute Conjecture~(A). Many efforts have been invested and we believe that will bear fruit one day. An account of the state of affairs in the Euler and Navier-Stokes equations around 2008 is due to P. Constantin \cite{Const08}. At the present moment  we are entertained by a number of trials and false proofs (some of them quite well published). There are excellent general texts on Navier-Stokes, such as \cite{Galdi} and~\cite{Temam}. Two very recent texts are \cite{RoRoSa16} and~\cite{Ser15}.

\section{Elliptic Equations and the Calculus of Variations}

For reasons of selection and space, we will be quite brief on a subject in which Louis made so many contributions.
We mention first of all the article \cite{BN83} by Haim Brezis and Louis Nirenberg, which figures among the most widely read among the works of both authors. It deals with the existence of solutions of semilinear elliptic equations with critical exponent (once again!)
\[
\Delta u + f(x,u)+ u^{(n+2)/(n-2)}=0.
\]
Two further articles that had great impact are work in collaboration with Shmuel Agmon and Avron Douglis
\cite{ADN59}, year 1959, and \cite{ADN64}, year 1964. They are near-the-boundary estimates for solutions of elliptic equations that satisfy general boundary conditions. Behavior near the boundary of non-linear or degenerate PDE solutions, or in domains with non-smooth boundaries, is a really delicate issue. Indeed, it is a topic of permanent interest in our community, in theory and also because of its practical interest\footnote{Think about the behavior of fluids in domains with corners.}.

The article \cite{BNV94} with Henri Berestycki and S.\,R.\,S. Varadhan links the study of the first eigenvalue with the maximum principle, a subject that Louis enjoyed so much. In this context we find the famous article on the method of the ``moving planes'' of 1991 \cite{BN91} in collaboration with Henri Berestycki, which I consider a gem.

In the Calculus of Variations let us quote the article \cite{BN83} with Haim Brezis, about the difference between local minimizers in the spaces $H^1$ and~$C^1$. See also~\cite{BrN91}.

A topic of great interest for Louis was the study of geometric properties such as symmetry.
The articles \cite{GNN79, GNN81} with Vasilis Gidas and Wei-Ming Ni deal with the radial symmetry of certain positive solutions of nonlinear elliptic equations that is imposed by the equation and the shape of the domain.

We are sorry for such an unjustly short space given to great contributions.
%
\section{Other contributions}

We collect here brief comments on important results obtained by Louis and his collaborators
on various topics that would deserve a more extensive treatment, for which we apologize to the expert reader.

\subsection{Operator theory}

Nirenberg and Joseph J. Kohn\footnote{J. J. Kohn is a brilliant Princeton analyst, not to be confused with R. Kohn from Courant. J. J. Kohn speaks perfect Spanish with an Ecuadorian accent.} introduced the notion of a pseudo-differential operator that helped generate a huge amount of later work in the brilliant school of harmonic analysis. In a 1965 article, \cite{KN65}, they dealt with pseudo-differential operators with a complete and algebraic view. The operators in question act on the space of tempered distributions at $ \R^n $, and are estimated in terms of Fourier transform norms.
The importance of these results is that they take into account all the ``lower order terms'', difficult to deal with in previous articles. See also the volume \cite{Nir10} edited by Louis.

\subsection{Free boundary problems}

This is one of the favorite topics of this reviewer. In 1977 Louis published with David Kinderlehrer the article \cite{KinNir1} on the regularity of free boundary problems for elliptic equations, at the beginning of an era that was to witness great progress. To put it clearly, let us assume that $ u $ is a solution to the problem
\[
\Delta u\le f, \quad u\ge 0, \quad (\Delta u-f)\,u=0
\]
defined in a domain $D \subset \R^n$. Boundary data are also given at the fixed boundary $ \partial D$. These data are intended to determine not only $ u $ but also the positivity domain $\Omega = \{x\in D : u (x)> 0\}$, or still better the boundary of $\Omega$ that lies within $D$, called the \sl free boundary\rm:
\[
\Gamma(u) = \partial \Omega\cap D.
\]
This is properly called an \textit{obstacle problem}. To get a physical idea, we can imagine a membrane in space $\R^3$ of height $z = U (x, y)$ that is subject to boundary conditions $U = h \ge 0$ in $\partial D$ and must lie above a table (obstacle) of height $U_{\textup{obst}} (x, y) = 0$.

Often, we want to consider a nontrivial obstacle $\varphi$, usually a concave function as in the figure. This leads to an interesting equivalent formulation. If we put $u=U+\varphi$, we arrive at the problem
\[
\Delta u\le g, \quad u\ge \varphi, \quad (\Delta u-g)(U-\varphi) = 0,
\]
with driving term $g=f+\Delta \varphi$, and then we usually take $g=0$. In this formulation, $u$ is constrained to stay above the obstacle $u_{\textup{obst}}(x) = \varphi$.

In any case, in the ``free part'', $\{x\in \R^n: U(x)> 0 \} = \{x\in \R^n: u(x)> \varphi \}$, an elastic equation $\Delta U = f$ is satisfied, but a priori we do not know where that part could be located. It is therefore a problem that combines PDEs and Geometry (again!).

This problem was known to have a unique \sl solution pair\rm, $(u, \Gamma)$. The attentive reader will have observed that once $\Gamma$ is known, and with it $\Omega$, the PDE problem to find $u$ is rather elementary. Therefore, the difficulty lies in principle in the geometry. However, the solution to the puzzle was rather found in nonlinear analysis, \cite{KinSt80}, which also produces efficient numerical methods.

\begin{figure}
\centering
\includegraphics[width=0.45\textwidth]{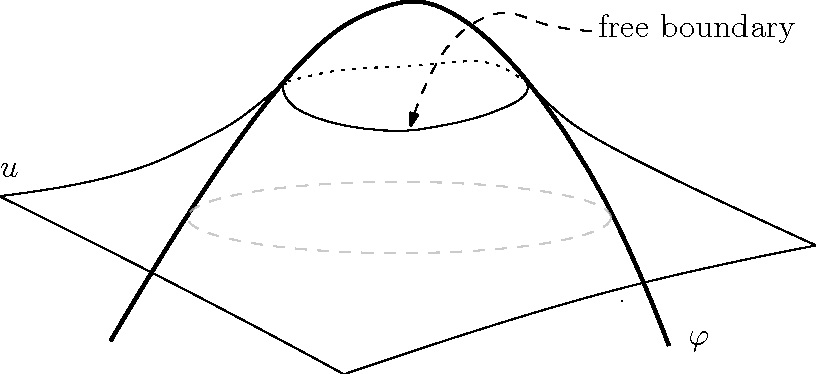}
\quad
\includegraphics[width=0.4\textwidth]{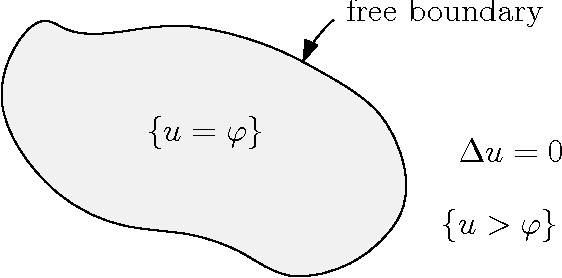}
\pseudocaption{Free boundaries and obstacles (pictures: X. Ros-Oton).}
\end{figure}

We then encounter a big theoretical problem: determining how regular is the set $\Gamma $, that we have found by abstract methods, and also determining how regular is $u$ near~$\Gamma$. Even the simplest question: ``is $\Gamma$ a surface?'' has to be answered. D. Kinderlehrer and L. Nirenberg gave local conditions on $f$ and assumed a certain initial regularity of $u$ to conclude that then $\Gamma$ is a very regular, even analytical, hyper-surface. The study of free boundaries extends to problems evolving in time, such as the very famous Stefan problem discussed by Louis in~\cite{KinNir2}. The 1980s were years of great progress in the mathematical understanding of free boundaries, with reference books such as \cite{Diaz} and~\cite{Fried}.

This is a field of very intense activity, both theoretical and applied, in which I have worked with great delight for decades. A required reference for in-depth study of the regularity of free boundaries is the book \cite{CaffSalsa} by L. Caffarelli and S. Salsa, see also A. Petrosyan et al.~\cite{PSU12}. A study of tumor growth modeling, seen as a free boundary problem, was done by B. Perthame et al.\ in \cite{PQV14}, it is just an example from a vast literature.

\subsection{Geometric Equations}

The article \cite{LN74} with Charles Loewner in 1974 deals with PDEs that are invariant under conformal or projective transformations. The reader will recall in this context the current relevance of PDEs linked to problems of Riemannian geometry, such as the Yamabe problem. We refer to the lengthy overview \cite{YYli} due to Yan Yan Li, Louis's doctoral student that has been for many years professor at Rutgers.

\subsection{Complex geometry}

The topic interested Louis a lot in his beginnings. A. Newlander and L. Nirenberg wrote in 1965 an article published in Annals of Mathematics \cite{NN57} on analytical coordinates in quasi-complex manifolds. The Newlander-Nirenberg Theorem states that any integrable quasi-complex structure is induced by a complex structure. Integrability is expressed through a list of differential conditions.

\medskip

\begin{center}
$\clubsuit$
\end{center}

\medskip

We put an end here to the mathematical journey, unfortunately unfair in many aspects due to the brevity of space and my ignorance in so many subjects. We hope that the extensive cited literature will serve as an indication of the profound influence of Louis Nirenberg and his world on the mathematicians and mathematics that have followed him. For the curious reader, there are excellent articles dealing with the work and life of Louis Nirenberg: a congress in his honor on the occasion of the 75th anniversary was organized by Alice Chang et al.\ and is collected in~\cite{ChLinYau}. He was interviewed by Allyn Jackson for the \textit{AMS Notices} in 2002, \cite{AJack}, and Simon Donaldson,
Fields Medal, wrote about him in the same journal in 2011, \cite{SD2011}. Yan-Yan Li's \cite {YYli} 2010 article focuses on the analysis of geometric problems. On the occasion of the Abel Prize, Xavier Cabré wrote a review in Catalan in \cite{XC15} and Tristan Rivi\`ere reviews his work in PDEs in \cite{Riv2015}.
A mathematical description of the influence of his ideas appeared in 2016 in \cite{Noti16} with contributions of a number of experts: X. Cabr\'e (symmetries of solutions), A. Chang (Gauss curvature problem), G. Seregin (Navier-Stokes problem), E. Carlen and A. Figalli (stability of the GNS inequality), M. T. Wang and S. T. Yau (Weyl problem and general relativity). Finally, the book \cite{Holden2019} presents the laureates of the Abel Prize in the period 2013--2017. In it Robert V. Kohn devotes to L. Nirenberg the article ``A few of Louis Nirenberg's many contributions to the theory of partial differential equations''.
By the way, there is a beautiful quotation from Abel as motto for the book: \textit{``Au reste il me para\^it que si l'on veut faire des progr\`es dans les math\'ematiques il faut \'etudier les ma\^itres et non pas les \'ecoliers''}.\footnote{In English: ``Finally, it appears to me that if one wants to make progress in mathematics, one should study the masters, not the students.'' Taken from the book.}

\medskip

\noindent {\bf Update.} We are happy to acknowledge the article ``A personal tribute to Louis Nirenberg'',
posted by Prof. Joel Spruck in the Arxiv repository in May 2021, \cite{Spruck2021}. As a person who met Louis Nirenberg in 1972 and  became a Courant Instructor, his detailed report on a selection of Louis's works is a very commendable reading. He concentrates on the work inspired by geometric problems beginning around 1974, especially the method of moving planes, and implicit fully nonlinear elliptic equations, and makes comments on Louis' personality.

\section{The quiet wise man and Spain}

I hope the reader will allow me to conclude this essay by some personal notes, that bear relation to Spain. My first memory of Louis Nirenberg sets us in Lisbon in the spring of 1982.\footnote{At the International Symposium in Homage to Prof.\ J. Sebati\~ao e Silva.} He was already famous and I was a novice in the art. In Lisbon I listened to one of his talks, which brought together the depth of the mathematics, the simplicity of the exposition and a grace to add some comment as timely as it was nice, characteristic features of Louis that delighted the public.

In the fall of that same year I set foot in the US, headed for the University of Minnesota,\footnote{This American university was very popular with young Spanish graduates and doctors for the excellence of its studies in Mathematics and Economics.} to work on free boundary problems with Don Aronson and with Luis Caffarelli, who was back from his visit to Courant Institute. Then I saw, through the group of great professor I had access to, that mathematical research offered a much better way of life. Among that group of friends I count Haim Brezis and Luis Caffarelli who have been my masters, Louis Nirenberg, Constantine Dafermos, Donald Aronson, Mike Crandall, Hans Weinberger,\dots\ I will never cease from thanking them for that vision.

\begin{figure}
\centering
\includegraphics[width=0.82\textwidth]{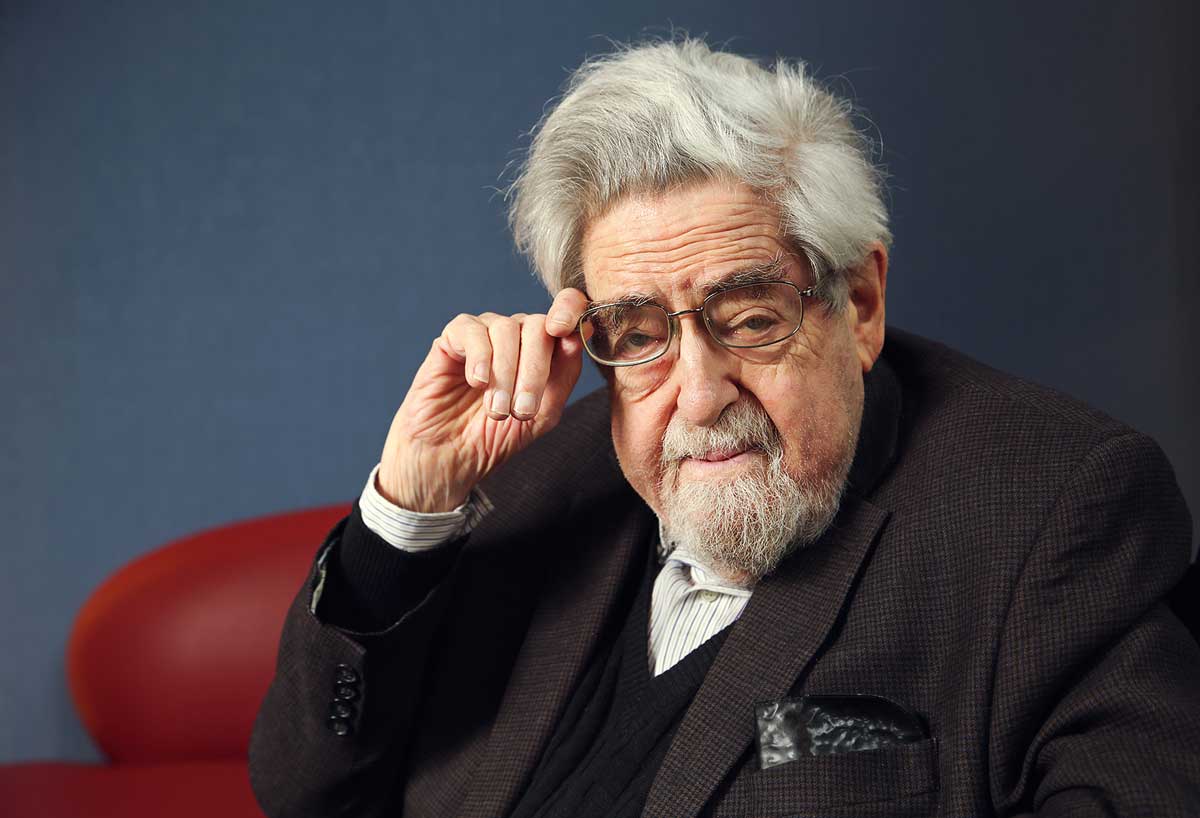}
\pseudocaption{Nirenberg in Barcelona in 2017 (photo: Jordi Play).}
\end{figure}

A few years later, I had the honor of participating in the organization of a summer course at the UIMP\footnote{Menéndez Pelayo International University, the course took place in 1987 at the Palacio de la Magdalena in Santander.} which included Louis as lecturer along with Don G. Aronson (Minnesota), Philippe Bénilan (Besançon), Luis A. Caffarelli (IAS Princeton) and Constantine Dafermos (Brown Univ.). These courses were inspired by Luis Caffarelli, close collaborator and friend of Louis, with the support of the Rector of the UIMP, Prof.\ Ernest Lluch,\footnote{Scholar of indelible memory, great protector of science and great conversationalist, he died tragically for being a good person at a very turbulent time.} and somehow they transmitted a certain spirit of mathematics that was being done around the Courant Institute. The course had a remarkable consequence. A young mathematician from Barcelona, Xavier Cabré,
a student in the course, went to the Courant Institute with Louis Nirenberg and thus began an international mathematical career, like the ones that so many young people crave today. His thesis, directed by Louis, dealt with ``Estimates for Solutions of Elliptic and Parabolic Equations'' (NYU, 1994). Following his stay in New York, he published with Luis Caffarelli the beautiful book \cite{CC95} on the so-called completely nonlinear elliptic equations. Xavier Cabré is now an ICREA Professor at the UPC in Barcelona. Louis Nirenberg visited Spain several times, specially Barcelona, and had many Spanish friends and admirers.

Although I did not become a collaborator of Louis, I had the opportunity of seeing him and talking to him on several occasions. I highlight a stay at the Courant Institute in the winter of 1996 where I could appreciate the day-to-day life of the ``quiet wise man'', or a congress in Argentina in 2009 when he was already very senior but loved life as the first day. The last event in which I saw him took place at Columbia University, New York, in May of last year (2019), in a congress in honor of Luis Caffarelli. He went to some talks in his wheelchair at 94 years old, and, with his proverbial good humor he told us that it was a bit difficult for him to follow the lectures!

Impressed by his personality, the young mathematician David Fernández and myself wrote a portrait of him in two entries in the blog ``The Republic of Mathematics'' that we edit in ``Investigación y Ciencia'' (Spanish partner of ``Sciencific American''). We called the essays ``Louis Nirenberg, the quiet wise man'' (I) and (II).\footnote{\url{https://www.investigacionyciencia.es/blogs/matematicas/75/posts}.} He was a teacher and master of science as those described by George Steiner in \cite{St2004}, where the relationship between teacher and pupil, master and disciple, is what matters. Louis had 46 doctoral students, many of them well-known mathematicians.\footnote{The first was Walter Littman (in 1956), whom I treated so much in Minnesota.} It was not his style to write long textbooks, he was the author of \cite{Nir74} and the recently published~\cite{Nir18}.

We will miss the teacher, master and senior friend who always looked gentle and kind, who loved Italy (\textit{il bel paese}), culture, good food and talking about movies and friends, and with whom mathematics was easy and exciting. Nirenberg lived in New York since 1949, in the Upper West Side, he was a perfect New Yorker and at the same time a citizen of the wide world. He worked until the end of his life, frequently visiting ``his'' Institute. Lucky soul, how I envy him, now and here the ``elders'' seem expendable for public utility.

I am proud to bear his name Louis = Luis, like Luis Caffarelli or Jacques Louis Lions or Luigi Ambrosio. He is already a great name in mathematics and it is an honor that carries the burden of working as Louis Nirenberg, only for the best and always in a good mood, and that is not easy.
Rest in eternal peace, beloved Master. In the Elysian fields you will have time to think about new functional inequalities, the beautiful functions that optimize them, and their surprising fruits. In our own small way, we also follow them, as in~\cite{DTGV}.


\section*{Acknowledgments and credits}

This is an extended English version of the article appeared in \textit{La Gaceta de la RSME}, the Spanish Royal Mathematical Society, \cite{VGac20}. I am deeply indebted to RSME for the invitation, continued interest, and technical support. In writing this essay I have used public domain documents from sites such as MathSciNet, Google Scholar, Wikipedia, Biographies of MacTutor History of Mathematics,\footnote{\url{http://mathshistory.st-andrews.ac.uk/Biographies/Nirenberg.html}.} The Mathematics Genealogy Project,\footnote{\url{https://www.genealogy.math.ndsu.nodak.edu/}.} the pages of the Abel Foundation\footnote{\url{https://www.abelprize.no/}.} and the Clay Mathematics Institute.\footnote{\url{https://www.claymath.org/}.}
Some paragraphs are adapted from my writings for the aforementioned blog, for the Oviedo newspaper \textit{La Nueva Espa\~na},\footnote{\url{https://www.lne.es/sociedad/2020/01/28/muere-louis-nirenberg-genio-matematicas/2590290.html}.}
or for the \textit{Boletín de la~RSME.}\footnote{\url{https://www.rsme.es/wp-content/uploads/2020/01/Boletin653.pdf}.}

The author acknowledges the support of a number of collaborators who have provided data, suggestions and corrections.




\end{document}